\documentclass[draft]{amsart}
\usepackage{amssymb}
\newcommand\Zset{\mathbb {Z}}
\newcommand\Qset{\mathbb {Q}}
\newcommand\Cset{\mathbb {C}}

\newcommand\tr{\mathrm{tr}}
\newcommand\A{{\mathcal A}}
\newcommand\U{{\mathcal U}}
\newcommand\B{{\mathcal B}}
\newcommand\la{l^2(\A)}
\newcommand\lge{l^2(G)}
\newcommand\vng{{\mathcal N}G}
\newcommand\ug{{\mathcal U}G}
\newcommand\NZ{{\mathcal N}\Zset}
\newcommand\UZ{{\mathcal U}\Zset}
\newcommand\ef{{\mathcal F}}
\newcommand\te{{\mathcal T}}
\newcommand\ce{{\mathcal C}}
\newcommand\T{\mathrm{{\bf T}}}
\newcommand\smallt{\mathrm{{\bf t}}}
\newcommand\bigP{\mathrm{{\bf P}}}
\newcommand\p{\mathrm{{\bf p}}}
\newcommand\unb{\mathrm{{\bf u}}}
\newcommand\bnd{\mathrm{{\bf b}}}
\newcommand\cl{\mathrm{cl}}
\newcommand\tor{\mathrm{Tor}}
\newcommand\homo{\mathrm{Hom}}

\newtheorem{thm}{Theorem}[section]
\newtheorem{lem}{Lemma}[section]

\newtheorem{prop}{Proposition}[section]
\newtheorem{defn}{Definition}[section]
\newtheorem{exmp}{Example}[section]

\begin{document}

\title{Torsion Theories for Algebras of Affiliated Operators
of Finite von Neumann Algebras
}

\author{Lia Va\v s}

\address{Department of Mathematics, Physics and Computer Science,
University of the Sciences in Philadelphia, 600 S. 43rd St.,
Philadelphia, PA 19104}

\email{l.vas@usip.edu}

\thanks{Part of the results are obtained during the time the author
was at the University of Maryland, College Park. The author was
supported by NSF grant DMS9971648 at that time.}

\subjclass[2000]{16W99, %None of the above, but in this section (16W is Rings and algebras with additional structure)
46L99, %None of the above, but in this section 46L is Selfadjoint operator algebras
16S90} %Maximal ring of quotients, torsion theories, radicals on module categories


\begin{abstract}
The dimension of any module over an algebra of affiliated
operators $\U$ of a finite von Neumann algebra $\A$ is defined
using a trace on $\A.$ All zero-dimensional $\U$-modules
constitute the torsion class of torsion theory $(\T,\bigP)$. We
show that every finitely generated $\U$-module splits as the
direct sum of torsion and torsion-free part. Moreover, we prove
that the theory $(\T,\bigP)$ coincides with the theory of bounded
and unbounded modules and also with the Lambek and Goldie torsion
theories. Lastly, we use the introduced torsion theories to give
the necessary and sufficient conditions for $\U$ to be semisimple.
\end{abstract}

\maketitle

\section{{\bf Introduction}}

A finite von Neumann algebra proves to be an interesting structure
both for operator theorists and for those working in geometry or
algebraic $K$-theory. One of the reasons is that a finite von
Neumann algebra $\A$ comes equipped with a normal and faithful
trace that enables us to define the dimension not just of a
finitely generated projective module over $\A$ but also of
arbitrary $\A$-module.

Moreover, $\A$ mimics the ring $\Zset$ in such a way that every
finitely generated module is a direct sum of a torsion and
torsion-free part. The dimension faithfully measures the
torsion-free part and vanishes on the torsion part. $\A$ has nice
ring-theoretic properties: it is semihereditary (i.e., every
finitely generated submodule of a projective module is projective)
and an Ore ring. The fact that $\A$ is an Ore ring allows us to
define the classical ring of quotients denoted $\U$. Besides this
algebraic definition, it turns out that, within the operator
theory, $\U$ can be defined as the algebra of affiliated
operators.

Using the dimension over $\A$, we can define the dimension over
$\U$ with the same properties as the dimension over $\A.$ As a
ring, $\U$ keeps all the properties of the ring $\A$ and possesses
some additional properties that $\A$ does not necessarily have. In
the analogy that $\A$ is like $\Zset,$ $\U$ plays the role of
$\Qset.$ In Section \ref{VNAandU}, we define a finite von Neumann
algebra $\A$, the dimension of $\A$-module and the algebra of
affiliated operators of $\A,$ and list some results on these
notions that we shall use further on.

Every finitely generated module over a finite von Neumann algebra
$\A$ is a direct sum of a torsion and a torsion-free module.
However, it turns out that there exists more than just one
suitable candidate when it comes to defining torsion and
torsion-free modules. To overcome this problem, the notion of a
torsion theory of a ring comes in as a good framework for the
better understanding of the structure of $\A$-modules. In Section
\ref{TT}, we define a torsion theory for any ring and some related
notions. We introduce some torsion theories for a finite von
Neumann algebra $\A$: Lambek, Goldie, classical torsion theory,
the torsion theory $(\T,\bigP)$ (studied also in \cite{Lu1},
\cite{Lu2}, \cite{Tele_prvenac}, \cite{Re} for finitely generated
modules and in \cite{Tele_Teza} for group von Neumann algebras) in
which a module is torsion if its $\A$-dimension is zero and,
finally, the torsion theory $(\bnd,\unb)$ of bounded and unbounded
modules.

In Section \ref{TTforU}, we study the torsion theories for $\U$.
Since the dimension of a $\U$-module can be defined via the
dimension over $\A$, we can define the torsion theory $(\T,\bigP)$
in the same way as for $\A.$ If $M$ is a finitely generated
$\U$-module, we show that the short exact sequence
\[0\rightarrow\T M \rightarrow M \rightarrow\bigP M\rightarrow0\]
splits just as for finitely generated $\A$-modules (Proposition
\ref{b splits in U}). Then we show (Theorem \ref{T=b for U}) that,
for $\U$, \[(\T, \bigP) = \mbox{ Lambek torsion theory = Goldie
torsion theory}= (\bnd, \unb).\] This indicates that, in contrast
to the situation with $\A,$ there is only one nontrivial torsion
theory of interest for $\U.$

Thus, one can work with $\U$ instead of $\A$ if one is not
interested in the information that gets lost by the transfer from
$\A$ to $\U$ (faithfully measured by the Novikov-Shubin invariant,
see \cite{Lu5}). For applications to topology, see section 4 in
\cite{Re} or chapter 8 in \cite{Lu5} for details.

The passage from $\A$ to $\U$ mimics in many ways the passage from
a principal ideal domain to its quotient field. However, although
$\U$ has many nice properties as a ring, it is not necessarily
semisimple. Any infinite group gives us the group von Neumann
algebra with algebra of affiliated operators that is not
semisimple (see Exercise 9.11 in \cite{Lu5}). In Section
\ref{semisimple}, we use the introduced torsion theories to give
necessary and sufficient conditions for $\U$ to be semisimple
(Theorem \ref{ug semisimple second set}).

\section{{\bf Finite von Neumann Algebras and the Algebras of
Affiliated operators}} \label{VNAandU}

Let $H$ be a Hilbert space and $\B(H)$ be the algebra of bounded
operators on $H$. The space $\B(H)$ is equipped with five
different topologies: norm, strong, ultrastrong, weak and
ultraweak . The statements that a $*$-closed unital subalgebra
$\A$ of ${\B}(H)$ is closed in weak, strong, ultraweak and
ultrastrong topologies are equivalent (for details see \cite{D}).

\begin{defn} A {\em von Neumann algebra} $\A$ is a $*$-closed unital
subalgebra of ${\B}(H)$ which is closed with respect to weak
(equivalently strong, ultraweak, ultrastrong) operator topology.
\end{defn}

A $*$-closed unital subalgebra $\A$ of ${\B}(H)$ is a von Neumann
algebra if and only if $\A= \A''$ where $\A'$ is the commutant of
$\A.$ The proof can be found in \cite{D}.

\begin{defn} A von Neumann algebra $\A$ is {\em finite} if
there is a $\Cset$-linear function $\tr_{\A}: \A\rightarrow\Cset$
such that
\begin{enumerate}
\item $\tr_{\A}(ab) =\tr_{\A}(ba).$
\item $\tr_{\A}(a^*a)\geq 0$. A $\Cset$-linear function on $\A$
that satisfies 1. and 2. is called a {\em trace}.
\item $\tr_{\A}$ is {\em normal}: it is continuous with respect to
ultraweak topology.
\item $\tr_{\A}$ is {\em faithful}: $\tr_{\A}(a)=0$ for some $a\geq 0$
(i.e. $a = bb^*$ for some $b\in\A$) implies $a=0$.
\end{enumerate}
 \label{finite vna}
\end{defn}

A trace on a finite von Neumann algebra is not unique. A trace
function extends to matrices over $\A$ in a natural way: the trace
of a matrix is the sum of the traces of the elements on the main
diagonal.

\begin{exmp}
Let $G$ be a (discrete) group. The group ring $\Cset G$ is a
pre-Hilbert space with an inner product: $ \langle \;\sum_{g\in G}
a_g g, \sum_{h\in G} b_h h\;\rangle = \sum_{g\in
G}a_g\overline{b_g}.$

Let $\lge$ be the Hilbert space completion of $\Cset G.$ Then
$\lge$ is the set of square summable complex valued functions over
the group $G$.

{\em The group von Neumann algebra} $\vng$ is the space of
$G$-equivariant bounded operators from $\lge$ to itself:
\[\vng = \{\; f\in {\B}(\lge)\;|\;f(gx)=gf(x)\mbox{ for
all }g\in G\mbox{ and }x\in\lge\;\}.\]

$\Cset G$ embeds into ${\B}(\lge)$ by right regular
representations. $\vng$ is a von Neumann algebra for $H = \lge$
since it is the weak closure of $\Cset G$ in ${\B}(\lge)$ so it is
a $*$-closed subalgebra of ${\B}(\lge)$ which is weakly closed
(see Example 9.7 in \cite{Lu5} for details). $\vng$ is finite as a
von Neumann algebra since it has a normal, faithful trace
$\tr_{\A}(f) = \langle f(1), 1 \rangle.$
\end{exmp}

The trace provides us with a way of defining a convenient notion
of the dimension of any module.

\begin{defn}
If $P$ is a finitely generated projective $\A$-module, there exist
$n$ and $f:\A ^n\rightarrow\A^n$ such that $f=f^2=f^*$ and the
image of $f$ is $P.$ Then, the {\em dimension} of $P$ is \[
\dim_{\A}(P)=\tr_{\A}(f)\in[0,\infty).\] Here the map $f^*$ is
defined by transposing and applying $*$ to each entry of the
matrix corresponding to $f.$

If $M$ is any $\A$-module, the {\em dimension} $\dim_{\A}(M)$ is
defined as \[\dim_{\A}(M)=\sup \{ \dim_{\A}(P)\; |\; P \mbox{ fin.
gen. projective submodule of }M\}\in[0, \infty]. \]
\label{dimension of any vng module}
\end{defn}

The dimension of a finitely generated projective module $P$ does
not depend on the choice of $f$ and $n$ from the definition above
and depend only on the isomorphism class of $P.$ For more details,
see comments following Assumption 6.2 on page 238 of \cite{Lu5} or
remarks following the Definition 1.6 in \cite{Lu2}.

The dimension of arbitrary module is also well defined by Theorem
0.6 from \cite{Lu2} or, equivalently Theorem 6.5 and 6.7 from
\cite{Lu5}.

The dimension has the following properties.
\begin{prop}
\begin{enumerate}
\item If $\;0\rightarrow M_0\rightarrow M_1\rightarrow M_2\rightarrow
0$ is a short exact sequence of $\A$-modules, then $
\dim_{\A}(M_1)= \dim_{\A}(M_0)+\dim_{\A}(M_2).$
\item If $M = \bigoplus_{i\in I} M_i,$ then $\dim_{\A}(M) =
\sum_{i \in I}\dim_{\A}(M_i).$
\item If $M = \bigcup_{i \in I}M_i$ is a directed union of
submodules, then $\dim_{\A}(M) = \sup\{\;\dim_{\A}(M_i)\; |\; i\in
I\;\}.$
\item If $M$ is finitely generated projective, then
$\;\dim_{\A}(M)=0$ iff $M=0.$
\end{enumerate}
\label{propofdim}
\end{prop}

The proof of this proposition can be found in \cite{Lu2} or
\cite{Lu5}.

As a ring, a finite von Neumann algebra $\A$ is {\em
semihereditary} (i.e., every finitely generated submodule of a
projective module is projective or, equivalently, every finitely
generated ideal is projective). This follows from two facts.
First, every von Neumann algebra is an $AW^*$-algebra and, hence,
a Rickart $C^*$-algebra (see Chapter 1.4 in \cite{Be2}). Second, a
$C^*$-algebra is semihereditary as a ring if and only if it is
Rickart (see Corollary 3.7 in \cite{AG}). The fact that $\A$ is
Rickart also gives us that $\A$ is  {\em nonsingular} (see 7.6 (8)
and 7.48 in \cite{Lam}).

Note also that every statement about left ideals over $\A$ can be
converted to an analogous statement about right ideals. This is
the case because $\A$ is a ring with involution (which gives a
bijection between the lattices of left and right ideals and which
maps a left ideal generated by a projection to a right ideal
generated by the same projection).

\subsection{The Algebra of Affiliated Operators}

A finite von Neumann algebra is a pre-Hilbert space with inner
product $\langle a, b \rangle = \tr_{\A}(ab^*).$ Let $\la$ be the
Hilbert space completion of $\A.$ Note that in the group case
$l^2(\vng)$ is isomorphic to $\lge$ since they are both Hilbert
space completions of $\vng$ (see section 9.1.4 in \cite{Lu5} for
details). Also, a finite von Neumann algebra $\A$ can be
identified with the set of $\A$-equivariant bounded operators on
$\la,$ ${\B}(\la)^{\A},$ using the right regular representations.
This justifies the definition of $\vng$ as $G$-equivariant
operators in ${\B}(\lge)$ since $ {\B}(l^2(\vng))^{\vng}=
{\B}(\lge)^{\vng} = {\B}(\lge)^{G} = \vng.$

\begin{defn} Let $a$ be a linear map $a: \mathrm{dom}\; a$
$\rightarrow$ $\la$, $\mathrm{dom}\; a \subseteq\la.$ We say that
$a$ is {\em affiliated to $\A$} if
\begin{itemize}
\item[i)] $a$ is densely defined (the domain $\mathrm{dom}\; a$
is a dense subset of $\la);$
\item[ii)] $a$ is closed (the graph of $a$ is closed in
$\la \oplus \la);$
\item[iii)]  $ba = ab$ for every $b$ in the commutant of $\A.$
\end{itemize}
Let $\U=\U(\A)$ denote the {\em algebra of operators affiliated
to} $\A$.
\end{defn}

\begin{prop}
Let $\A$ be a finite von Neumann algebra and $\U$ its algebra of
affiliated operators.
\begin{enumerate}
\item $\A$ is an Ore ring. \item $\U$ is equal to the classical
ring of quotients $Q_{\mathrm{cl}}(\A)$ of $\A.$ \item $\U$ is a
von Neumann regular, left and right self-injective ring. \item
$\U$ is the maximal ring of quotients $Q_{\mathrm{max}}(\A).$
\end{enumerate}
\label{ug is classical and maximal}
\end{prop}
The proof of 1. and 2. can be found in \cite{Re}. The proof of 3.
and 4. can be found in \cite{Be1}.

From this proposition it follows that the algebra $\U$ can be
defined using purely algebraic terms (ring of quotient, injective
envelope) on one hand and using just the language of operator
theory (affiliated operators) on the other.

The ring $\U$ has many nice properties that $\A$ is missing: it is
von Neumann regular and self-injective; and it keeps all the
properties that $\A$ has: it is semihereditary and nonsingular.

Further, $K_0(\A)$ and $K_0(\U)$ are isomorphic. Namely, Handelman
proved (Lemma 3.1 in \cite{Handel}) that for every finite Rickart
$C^*$-algebra $R$ such that every matrix algebra over $R$ is also
Rickart, the inclusion of $R$ into a certain regular ring $U(R)$
with the same lattice of projections as $R$ induces an isomorphism
$\mu: K_0(R)\rightarrow K_0(U(R)).$ By Theorem 3.4 in \cite{AG}, a
matrix algebra over a Rickart $C^*$-algebra is a Rickart
$C^*$-algebra. Thus, $K_0(R)$ is isomorphic to $K_0(U(R))$ for
every finite Rickart $C^*$-algebra $R.$ If $\A$ is a finite von
Neumann algebra, the regular ring from Handelman's theorem can be
identified with the maximal ring of quotients
$Q_{\mathrm{max}}(\A)$ (see \cite{Be1}). This gives us that the
inclusion of a finite von Neumann algebra $\A$ in its algebra of
affiliated operators $\U$ induces the isomorphism \[\mu:
K_0(\A)\rightarrow K_0(\U).\] In \cite{Tele_Teza} it is shown that
the inverse of this isomorphism is induced by the map Proj($\U$)
$\rightarrow$ Proj($\A$) given by $[Q]\mapsto[Q\cap\A^n]$ for any
direct summand $Q$ of  $\U^n.$ Thus, the following holds.

\begin{thm}
There is an one-to-one correspondence between direct summands of
$\A$ and direct summands of $\U$ given by $I \mapsto$
$\U\otimes_{\A}I = E(I).$ The inverse map is given by $L\mapsto
L\cap \A.$ This correspondence induces an isomorphism of monoids
$\mu :\mathrm{Proj}(\A)$ $\rightarrow$ $\mathrm{Proj}(\U)$ and an
isomorphism \[\mu: K_0(\A)\rightarrow K_0(\U)\] given by
$[P]\mapsto[\U\otimes_{\A}P]$ with the inverse
$[Q]\mapsto[Q\cap\A^n]$ if $Q$ is a direct summand of $\U^n.$
\label{moja K 0 theorem}
\end{thm}

In Chapter 4 of \cite{Tele_Teza}, this theorem was proved for a
group von Neumann algebra and in Chapter 7 of \cite{Tele_Teza}, it
is shown that it holds for any finite von Neumann algebra as well.
In \cite{Tele_prvenac}, this result is contained in Theorem 5.2.

In \cite{Re}, the dimension of an $\U$-module is defined using the
dimension of an $\A$-module and the above isomorphism $\mu.$ The
{\em dimension} over $\U$ of a finitely generated projective
$\U$-module $M$ is defined as \[ \dim_{\U}(M) =
\dim_{\A}(\mu^{-1}(M)),
\] where $\dim_{\A}(\mu^{-1}(M))$ denotes the dimension over $\A$
of any module in the inverse image of the equivalence class $[M].$

Just as for the ring $\A,$ we can extend the definition of the
dimension to all modules. If $M$ is an $\U$-module, define the
{\em dimension} of $M$, $\dim_{\U}(M),$ as follows: \[\dim_{\U}(M)
= \sup \{\dim_{\U}(P)\;|\;P\mbox{ is a fin. gen proj. submodule of
}M\}.\] The dimension over $\U$ is well defined. For details, see
section 8.3 in \cite{Lu5} or section 3 in \cite{Re}.

In \cite{Re}, it is shown that the dimension over $\U$ has all the
properties that the dimension over $\A$ had, i.e. Proposition
\ref{propofdim} holds for $\dim_{\U}$ as well. In addition, in
\cite{Re} it is shown that the following holds:
\[\dim_{\U}(\U\otimes_{\A}N) = \dim_{\A}(N)\mbox{ for every
}\A\mbox{-module }N.\]

In \cite{Tele_drugi}, it is shown that $\U$ allows the definition
of another type of dimension. This dimension is analogous to the
cental-valued dimension over a finite von Neumann algebra
considered in \cite{Lu1}. For more details, see Section 4.2 of
\cite{Tele_drugi}.

\section{{\bf Torsion Theories}}
\label{TT}

To study different ways of defining the torsion and torsion-free
parts of modules over $\A$ or $\U,$ we first introduce the general
framework in which we shall be working --- the torsion theory.

\subsection{General Torsion Theories}

\begin{defn} Let $R$ be any ring. A {\em torsion theory} for
$R$ is a pair $\tau = (\te, \ef)$ of classes of $R$-modules such
that
\begin{itemize}
\item[i)] $\homo_R(T,F)=0,$ for all $T \in \te$ and $F \in \ef.$
\item[ii)] $\te$ and $\ef$ are maximal classes having property
$i).$
\end{itemize}
\end{defn}
The modules in $\te$ are called {\em $\tau$-torsion modules} (or
torsion modules for $\tau$) and the modules in $\ef$ are called
{\em $\tau$-torsion-free modules} (or torsion-free modules for
$\tau$).

If $\tau_1 = (\te_1, \ef_1)$ and $\tau_2 = (\te_2, \ef_2)$ are two
torsion theories, we say that $\tau_1$ is {\em smaller} than
$\tau_2$ and write $\tau_1\leq\tau_2$ if and only if
$\te_1\subseteq\te_2$. Equivalently, $\tau_1\leq\tau_2$ iff
$\ef_1\supseteq\ef_2.$

If $\ce$ is a class of $R$-modules, then the torsion theory {\em
generated} by $\ce$ is the smallest torsion theory $(\te, \ef)$
such that $\ce\subseteq\te.$ The torsion theory {\em cogenerated}
by $\ce$ is the largest torsion theory $(\te, \ef)$ such that
$\ce\subseteq\ef.$

\begin{prop}
\begin{enumerate}
\item If $(\te, \ef)$ is a torsion theory, then the class $\te$ is closed under quotients, direct sums and
extensions and the class $\ef$ is closed under taking submodules,
direct products and extensions.

\item  If $\ce$ is a class of $R$-modules closed under quotients,
direct sums and extensions, then it is a torsion class for a
torsion theory $(\ce, \ef)$ where $\ef =
\{\;F\;|\;\homo_R(C,F)=0,\mbox{ for all }C\in\ce\;\}.$

Dually, if $\ce$ is a class of $R$-modules closed under
submodules, direct products and extensions, then it is a
torsion-free class for a torsion theory $(\te, \ce)$ where $\te =
\{\;T\;|\;\homo_R(T,C)=0,\mbox{ for all }C\in\ce\;\}.$
\item Two classes of $R$-modules $\te$ and $\ef$ constitute a
torsion theory iff
\begin{itemize}
\item[i)] $\te\cap\ef = \{0\},$
\item[ii)] $\te$ is closed under quotients,
\item[iii)] $\ef$ is closed under submodules and
\item[iv)] For every module $M$ there exists submodule $N$ such that
$N\in\te$ and $M/N\in\ef.$
\end{itemize}
\end{enumerate}
\end{prop}
The proof of this proposition is straightforward by the definition
of a torsion theory. The details can be found in \cite{Bland}.

From this proposition it follows that every module $M$ has a
largest submodule which belongs to $\te$. We call it the {\em
torsion submodule} of $M$ and denote it $\te M$. The quotient
$M/\te M$ is called the {\em torsion-free quotient} and we denote
it $\ef M.$

A torsion theory $\tau = (\te, \ef)$ is {\em hereditary} if the
class $\te$ is closed under taking submodules. A torsion theory is
hereditary if and only if the torsion-free class is closed under
formation of injective envelopes. Also, a torsion theory
cogenerated by a class of injective modules is hereditary and,
conversely, every hereditary torsion theory is cogenerated by some
class of injective modules. The details can be found in
\cite{Golan}.

A torsion theory enables us to define the closure of a submodule
in a module.
\begin{defn}
If $M$ is an $R$-module and $K$ a submodule of $M,$ let us define
the {\em closure} $\cl_{\te}^M(K)$ of $K$ in $M$ with respect to
the torsion theory $(\te, \ef)$ by \[\cl_{\te}^M(K) =
\pi^{-1}(\te(M/K))\mbox{ where } \pi\mbox{ is the natural
projection }M\twoheadrightarrow M/K.\]
\end{defn}

If it is clear in which module we are closing the submodule $K,$
we suppress the superscript $M$ from $\cl_{\te}^M(K)$ and write
just $\cl_{\te}(K)$. If $K$ is equal to its closure in $M,$ we say
that $K$ is {\em closed} submodule of $M$.

The closure has the following properties.
\begin{prop}
Let $(\te, \ef)$ be a torsion theory on $R$, let $M$ and $N$ be
$R$-modules and $K$ and $L$ submodules of $M$. Then
\begin{enumerate}
\item $\te M = \cl_{\te}(0).$
\item $\te(M/K) = \cl_{\te}(K)/K$ and
$\ef(M/K) \cong M/\cl_{\te}(K).$
\item If $K \subset L,$ then $\cl_{\te}(K)\subseteq\cl_{\te}(L).$
\item $K\subset\cl_{\te}(K)$ and
$\cl_{\te}(\cl_{\te}(K)) = \cl_{\te}(K).$
\item $\cl_{\te}(K)$ is the smallest closed submodule of
$M$ containing $K.$
\item If $(\te, \ef)$ is hereditary, then
$\cl_{\te}^K(K\cap L) = K \cap \cl^M_{\te}(L).$ If $(\te, \ef)$ is
not hereditary, just $\subseteq$ holds in general.
\item If $(\te_1, \ef_1)$ and $(\te_2, \ef_2)$ are two torsion
theories, then $(\te_1, \ef_1)\leq(\te_2, \ef_2)$ if and only if
$\cl_{\te_1}(K)\subseteq\cl_{\te_2}(K)$ for all $K.$
\end{enumerate}
\label{properties of closure}
\end{prop}

The proof follows directly from the definition of the closure.

\subsection{Examples}
\label{Examples}

\begin{enumerate}

\item {\em The trivial torsion theory} on Mod$_R$ is the torsion
theory $(0, \mathrm{Mod}_R).$

\item {\em The improper torsion theory} on Mod$_R$ is the torsion
theory $(\mathrm{Mod}_R, 0).$

\item The torsion theory cogenerated by the injective envelope $E(R)$
of $R$ is called the {\em Lambek torsion theory}. We denote it
$\tau_L.$ Since it is cogenerated by an injective module, it is
hereditary.

If the ring $R$ is torsion-free in a torsion theory $\tau$, we say
that $\tau$ is {\em faithful}. $\tau_L$ is faithful and it is the
largest hereditary faithful torsion theory.

\item The class of nonsingular modules over a ring $R$ is closed
under submodules, extensions, products and injective envelopes.
Thus, the class of all nonsingular modules is a torsion-free class
of a hereditary torsion theory. This theory is called the {\em
Goldie torsion theory} $\tau_G.$

The Lambek theory is smaller than the Goldie theory. This is the
case since $\tau_G$ is larger than any hereditary torsion theory
(see \cite{Bland}). Moreover, $\tau_L = \tau_G$ if and only if $R$
is a nonsingular ring (i.e. $\tau_G$ is faithful). Recall that a
finite von Neumann algebra is a nonsingular ring.

\item
If $R$ is an Ore ring with the set of regular elements $T$ (i.e.,
$Tr \cap Rt \neq 0,$ for every $t \in T$ and $r\in R$), we can
define a hereditary torsion theory by the condition that an
$R$-module $M$ is a torsion module iff for every $m\in M$, there
is a nonzero $t\in T$ such that $tm =0.$ This torsion theory is
called the {\em classical torsion theory of an Ore ring}.

This theory is faithful and so it is contained in $\tau_L$.

\item Let $R$ be a subring of a ring $S$. Let us look at a
collection $\te$ of $R$-modules $M$ such that $S\otimes_R M = 0.$
This collection is closed under quotients, extensions and direct
sums. Moreover, if $S$ is flat as an $R$-module, then $\te$ is
closed under submodules and, hence, defines a hereditary torsion
theory. In this case, denote the torsion theory by $\tau_S.$

From the definition of $\tau_S$ it follows that the torsion
submodule of a module $M$ in $\tau_S$ is the kernel of the natural
map $M\rightarrow S \otimes_R M,$ i.e. Tor$^R_1(S/R, M).$ Thus,
all flat modules are $\tau_S$-torsion-free. Since $R$ is flat,
$\tau_S$ is faithful, so $\tau_S\leq \tau_L.$

If a ring $R$ is Ore, then the classical ring of quotients
$Q_{\mathrm{cl}}^l(R)$ is a flat $R$-module and the set $\{\;m\in
M\;|\; rm = 0,\mbox{ for some nonzero-divisor
 }r\in R\;\}$ is equal to the kernel $\ker
(M\rightarrow Q_{\mathrm{cl}}^l(R) \otimes_R M).$ Hence the
torsion theory $\tau_{Q_{\mathrm{cl}}^l(R)}$ coincides with the
classical torsion theory of $R$ in this case.

Since $\U=Q_{\mathrm{cl}}(\A)$ (see Proposition \ref{ug is
classical and maximal}), $\U$ is a flat $\A$-module and
$\tau_{\U}$ is the classical torsion theory of $\A.$

\item All the torsion theories we introduced so far are hereditary.
Let us introduce a torsion theory that is not necessarily
hereditary. Let $(\bnd, \unb)$ be the torsion theory cogenerated
by the ring $R$ (thus this is the largest torsion theory in which
$R$ is torsion-free). We call a module in $\bnd$ a {\em bounded
module} and a module in $\unb$ an {\em unbounded module}.

Since $(\bnd, \unb)$ is cogenerated by $R$, the closure of a
submodule $K$ of an $R$-module $M$ is $\;\cl_{\bnd}^M(K) = \{\;x
\in M\; |\; f(x)=0,$  for every $f\in\homo_{R}(M,R)$ such that
$K\subseteq \ker f\;\}.$

The Lambek torsion theory $\tau_L$ is contained in the torsion
theory $(\bnd, \unb)$ because $R$ is $\tau_L$-torsion-free. There
is another interesting relation between $\tau_L$ and $(\bnd,
\unb)$ torsion theory. Namely,
\begin{center}
$M$ is a $\tau_L$-torsion if and only if every submodule of $M$ is
bounded.
\end{center}
This is a direct corollary of the fact that $\homo_R(M, E(R)) = 0$
if and only if $\homo_R(N, R)=0,$ for all submodules $N$ of $M,$
which is an exercise in \cite{Cohn1}.  Also, it is easy to show
that $(\bnd, \unb)$ = $\tau_L$ if and only if $(\bnd, \unb)$ is
hereditary.
\end{enumerate}

To summarize, for any ring $R$ we have the following relationship
for the torsion theories:
\begin{center} Trivial $\leq$ $\tau_L$ $\leq$ $\tau_G$
$\leq$ $(\bnd, \unb)$ $\leq$ Improper.
\end{center}
If $R$ is an Ore nonsingular ring, then
\begin{center}
Trivial $\leq$ Classical = $\tau_{Q_{\mathrm{cl}}(R)}$ $\leq$
$\tau_L$ $=$ $\tau_G$ $\leq$ $(\bnd, \unb)$ $\leq$ Improper.
\end{center}
The last is the situation for our finite von Neumann algebra $\A$
as well as its algebra of affiliated operators $\U$. In the
following, we shall examine the situation for $\A$ and $\U$ in
more details.

\subsection{Torsion Theories for Finite von Neumann Algebras}

Let us introduce some theories for finite von Neumann algebras and
compare them with the torsion theories from previous chapter.

\begin{enumerate}
\item We can define a hereditary torsion theory using the dimension
of an $\A$-module. For an $\A$-module $M,$ define $\T M$ as the
submodule generated by all submodules of $M$ of the dimension
equal to zero. It is zero-dimensional by property (3) of
Proposition \ref{propofdim}. So, $\T M$ is the largest submodule
of $M$ of dimension zero. Let us denote the quotient $M/ \T M$ by
$\bigP M$.

The class $ \T = \{M \in \mathrm{Mod}_{\A} | M = \T M\}$ is closed
under submodules, quotients, extensions and direct sums
(Proposition \ref{propofdim}). Thus, $\T$ defines a hereditary
torsion theory with torsion-free class equal to $\bigP = \{M \in
\mathrm{Mod}_{\A} | M = \bigP M\}.$

From the definition of this torsion theory it follows that
$\cl_{\T}(K)$ is the largest submodule of $M$ with the same
dimension as $K$ for every submodule $K$ of an module $M$. Also,
since $\A$ is semihereditary and a nontrivial finitely generated
projective module has nontrivial dimension, $\A$ is in $\bigP$ and
so the torsion theory $(\T, \bigP)$ is faithful.

\item The second torsion theory of interest is $(\bnd, \unb),$ the
largest torsion theory in which the ring is torsion-free. Since
$\A$ is torsion-free in $(\T, \bigP),$ we have that $(\T, \bigP)
\leq (\bnd, \unb).$

In \cite{Lu2} it is shown that $\T M=\bnd M$ for a finitely
generated $\A$-module $M,$ that $\bigP M$ is a finitely generated
projective module and that $M=\bigP M\oplus \T M.$ The proof can
also be found in \cite{Lu5}.

\item Let $(\smallt, \p)$ denote the classical torsion theory of
$\A.$ Since $\U=Q_{\mathrm{cl}}(\A),$ \[\smallt M = \ker (M
\rightarrow \U \otimes_{\A}M)= \tor_1^{\A}(\U/\A, M)\] for any
$\A$-module $M$ (see Examples (5) and (6) in Subsection
\ref{Examples}).

Let $\p M$ denote the torsion-free quotient $M/ \smallt M.$ From
example (6), it follows that all flat modules are torsion-free. In
\cite{Turnidge}, the torsion theory from example (6) is studied.
Since $\A$ is semihereditary and
$\U=Q_{\mathrm{cl}}(\A)=Q_{\mathrm{max}}(\A)$ is von Neumann
regular and $\A$-flat, from Turnidge's results in \cite{Turnidge},
it follows that the converse holds as well: a torsion-free module
is flat. Hence, an $\A$-module $M$ is flat if and only if $M$ is
in $\p.$
\end{enumerate}

Various torsion theories for $\A$ are ordered as follows:
\begin{center}
Trivial $\leq$ Classical = $(\smallt,\p)$ $\leq$ $\tau_{L}$ $=$
$\tau_{G}$ = $(\T, \bigP)$ $\leq$ $(\bnd, \unb)$ $\leq$ Improper.
\end{center}

The proof of $\tau_{L}$ $=$ $\tau_{G}$ = $(\T, \bigP)$ can be
found in Chapter 4 of \cite{Tele_Teza} for the case of group von
Neumann algebras. The proof for the more general case of finite
von Neumann algebras is the same as for the group von Neumann
algebras (see Chapter 7 of \cite{Tele_Teza}). Alternatively,
Proposition 4.2 in \cite{Tele_prvenac} contains this result. It is
interesting to note that this proposition shows that the torsion
theory $(\T, \bigP)$, defined via a normal and faithful trace
$\tr_{\A},$ is not dependent on the choice of such trace since
$(\T, \bigP)$ coincides with the Lambek and Goldie theories.

The inequality $(\smallt, \p)\leq \tau_L$ holds since $\A$ injects
in $\U \otimes_{\A}\A = \U,$ so $\A$ is torsion-free in $(\smallt,
\p)$ and $\tau_L$ is the largest hereditary theory in which $\A$
is torsion-free.

All of the above inequalities can be strict. For details, see
\cite{Lu5} or \cite{Tele_Teza}.

If $M$ is an $\A$-module, there is a filtration:
\[\underbrace{0\subseteq\smallt}_{\smallt
M}\underbrace{M\subseteq\T }_{\T\p M} \underbrace{M\subseteq
M.}_{\bigP M}\] This follows from the fact that the quotient $\T
M/\smallt M = \p\T M$ is isomorphic to the module $\T\p M.$ For
details see Proposition 4.3 and comments following it in
\cite{Tele_prvenac}.

\section{{\bf Torsion Theories for the Algebra of Affiliated
Operators}} \label{TTforU}

Let us turn to the torsion theories of the algebra of affiliated
operators $\U$ of a finite von Neumann algebra $\A.$

Since we have defined the dimension over $\U$ and it satisfies all
the properties given in Proposition \ref{propofdim}, we can define
the hereditary torsion theory $(\T, \bigP)$ for $\U$ in the same
way as we did for $\A:$ the torsion submodule $\T M$ of a
$\U$-module $M$ is the greatest submodule of $M$ with of dimension
zero. $\bigP M$ is the quotient $M/ \T M.$ The class of all
zero-dimensional modules $\T$ is closed under quotients,
submodules, extensions and direct sums by Proposition
\ref{propofdim}. Hence, $(\T, \bigP)$ is a hereditary torsion
theory over $\U.$ $(\T, \bigP)$ coincides with the torsion theory
defined via the dimension considered in \cite{Tele_drugi}. For
more details, see Corollary 24 and the two paragraphs following it
in \cite{Tele_drugi}.

The second theory of interest is $(\bnd, \unb),$ the torsion
theory cogenerated by the ring itself. The Lambek torsion theory
$\tau_L$ is cogenerated by the injective envelope of the ring, but
$\U$ is a self-injective ring, hence $\tau_L = (\bnd, \unb).$
Further, since $\U$ is also a nonsingular ring, $\tau_L = \tau_G$.

$\U$ has no finitely generated submodules of dimension zero
because $\U$ is semihereditary and the dimension of a projective
module is zero only if the module is trivial. Since the dimension
of a module is the supremum of the dimensions of its finitely
generated submodules, $\U$ has no nontrivial submodules of
dimension zero. Thus $\U$ is in $\bigP$ so $(\T, \bigP)$ is
faithful. This yields \[(\T, \bigP)\leq \tau_L = \tau_G = (\bnd,
\unb).\] We will show that $(\T, \bigP) = (\bnd, \unb).$ The proof
consists of three steps. Lemma \ref{T is b for fgp Umod} tells us
that $\cl_{\T}=\cl_{\bnd}$ on submodules of a finitely generated
projective module. Proposition \ref{b splits in U} tells us that
$\cl_{\T}=\cl_{\bnd}$ on submodules of a finitely generated
module. Proposition \ref{b splits in U} will also tell us that a
finitely generated $\U$-module has the same property as a finitely
generated $\A$-module: it is the direct sum of its $\T$-submodule
and $\bigP$-quotient. Theorem \ref{T=b for U} will then tell us
that $\T=\bnd.$

Let $L_{\mathrm{fg}}(\U^n)$ denote the lattice of finitely
generated submodules of $\U^n$. Since $\U$ is von Neumann regular,
this lattice coincides with the lattice of direct summands of
$\U^n$. In \cite{Re} it is shown that this is a complete lattice
in which the supremum and infimum of two direct summands are their
sum and intersection, respectively. Note that the intersection of
two finitely generated $\U$-modules is finitely generated since
$\U$ is a coherent ring.

\begin{lem}
Let $P$ be a finitely projective $\U$-module, and $K$ a submodule
of $P$.
\[
\begin{array}{cccl} \cl_{\T}(K) & = &
\bigcap & \{Q\subseteq P\;|\; Q\mbox{ is finitely generated and
}K\subseteq Q \} \\
 & = & \inf & \{Q\subseteq P\;|\; Q\mbox{ is finitely generated  and
}K\subseteq Q \}\\
 & = & \sup & \{Q\subseteq P\;|\; Q\mbox{ is finitely generated  and
}Q\subseteq K \}\; = \;\cl_{\bnd}(K)
\end{array}
\]
$\cl_{\T}(K)$ is finitely generated and projective, and
$\cl_{\T}(K)$ is a direct summand of $P.$ \label{T is b for fgp
Umod}
\end{lem}
Note that the infimum and supremum  in the lemma denote the
operations in the lattice $L_{\mathrm{fg}}(\U^n)$ for $P$ a direct
summand of $\U^n.$ Since this lattice is complete, these two
modules are finitely generated and, hence, projective. The fact
that $\cl_{\T}(K)$ is finitely generated projective will follow
from the equality with these two modules.

\begin{proof}
Let $I$ ($I$ for infimum) denote the module $\inf\{Q\subseteq
P\;|\; Q$ is finitely generated and $K\subseteq Q \},$ $S$ ($S$
for supremum) denote $\sup\{Q\subseteq P\;|\; Q$ is finitely
generated and $Q\subseteq K \}$ and $Int$ ($Int$ for intersection)
denote the module $\bigcap\{Q\subseteq P\;|\; Q$ is finitely
generated and $K\subseteq Q \}.$ The proof proceeds in five steps:
\begin{enumerate}
\item $S = Int$;
\item $I = Int$;
\item $S\subseteq\cl_{\T}(K)$;
\item $\cl_{\T}(K)\subseteq\cl_{\bnd}(K)$;
\item $\cl_{\bnd}(K)= S.$
\end{enumerate}

(1) and (3) are proven in \cite{Re}.

(2) $Int$ is finitely generated projective by 1. (since $S$ is).
So, $Int$ is the largest finitely generated projective module that
is contained in all the modules $Q\subseteq P$ such that $Q$ is
finitely generated and $K\subseteq Q.$ But that means that $Int$
is the infimum of the set $\{Q\subseteq P\;|\; Q$ is finitely
generated and $K\subseteq Q \}.$ So, $I =Int.$

(4) $\cl_{\T}(K)\subseteq\cl_{\bnd}(K)$ follows since $(\T,
\bigP)\leq(\bnd, \unb).$

(5) $S\subseteq\cl_{\bnd}(K)$ by (3) and (4). We shall show the
equality by showing that $\cl_{\bnd}(K)/S$ is trivial. Note that
$\cl_{\bnd}(K)$ is equal to the intersection of the submodules
$\ker f$ where $f\in\homo_{\U}(P, \U)$ is such that
$K\subseteq\ker f$ (by the definition of the torsion theory
$(\bnd, \unb)$). The image of such a map $f$ is finitely generated
(since $P$ is) and projective (as a finitely generated submodule
of $\U$). But then $0\rightarrow\ker f\rightarrow P\rightarrow
\mathrm{im} f \rightarrow 0$ splits and so $\ker f$ is also
finitely generated projective. Since the lattice of finitely
generated submodules of $P$ is complemented (and the infimum is
just the intersection) $\cl_{\bnd}(K)$ is finitely generated
projective as well.

Since both $\cl_{\bnd}(K)$ and $S$ are finitely generated
projective, $\cl_{\bnd}(K)/S$ is finitely presented. All modules
over a von Neumann regular ring are flat and all finitely
presented flat modules are projective (Theorem 4.21, Theorem 4.30
in \cite{Lam}). Thus, a finitely presented module over a von
Neumann regular ring is finitely generated projective. So
$\cl_{\bnd}(K)/S$ is projective.

Since $\cl_{\bnd}(K)/S = \homo_{\U}(\U, \cl_{\bnd}(K)/S)$ to show
$\cl_{\bnd}(K)/S=0$ it is sufficient to show $\homo_{\U}(\U,
\cl_{\bnd}(K)/S)=0.$ But in every von Neumann regular ring $R,$
for two projective modules $P$ and $Q$ the following holds:
\begin{center} $\homo_R(P, Q)=0$ iff $\homo_R(Q, P)=0$
\end{center}
(this fact can be found in \cite{Go}). So, to show $\homo_{\U}(\U,
\cl_{\bnd}(K)/S)=0,$ it is sufficient to show
$\homo_{\U}(\cl_{\bnd}(K)/S, \U)=0.$ Let $\widetilde{f}$ be in
$\homo_{\U}(\cl_{\bnd}(K)/S, \U).$ It uniquely determines a map
$f: \cl_{\bnd}(K)\rightarrow\U$. Since $\U$ is self-injective, the
map $\homo_{\U}(P, \U)\rightarrow\homo_{\U}(\cl_{\bnd}(K), \U)$ is
onto. So, we can extend $f$ to $\overline{f}$ in $\homo_{\U}(P,
\U)$. Since $K\subseteq \ker f,$ $K\subseteq\ker \overline{f}$ and
so $\cl_{\bnd}(K)\subseteq\ker \overline{f}.$ But that means that
$\overline{f}|_{\cl_{\bnd}(K)}= f :\cl_{\bnd}(K)\rightarrow\U$ is
0, and so $\widetilde{f}=0$ as well. Hence,
$\homo_{\U}(\cl_{\bnd}(K)/S, \U)=0,$ which finishes the proof of
(5).

Since $\U$ is a self-injective ring, a finitely generated
projective module is injective. Thus, $\cl_{\T}(K)$ is a direct
summand of $P$ since it is finitely generated projective and a
submodule of $P.$
\end{proof}

The next proposition will tell us that $\cl_{\T}(K) =
\cl_{\bnd}(K)$ for every submodule $K$ of a finitely generated
$\U$-module $P.$

\begin{prop}
If $M$ is a finitely generated $\U$-module and $K$ a submodule of
$M,$ then
\begin{enumerate}
\item $\cl_{\bnd}(K)$ is a direct summand of $M$ and
$M/\cl_{\bnd}(K)$ is finitely generated projective.
\item $\dim_{\U}(K) = \dim_{\U}(\cl_{\bnd}(K)).$
\item $M = \bnd M\oplus \unb M$ and $\dim_{\U}(\bnd M)=0.$
\item $\T M=\bnd M$ so $M = \T M\oplus \bigP M$ and
$\bigP M = \unb M$ is a finitely generated projective module.
\end{enumerate}
\label{b splits in U}
\end{prop}
\begin{proof}
(1) Choose a finitely generated projective module $P$ and a
surjection $f: P\rightarrow M.$ By the previous lemma we know that
the $\T$-closure of a submodule in $P$ is the same as
$\bnd$-closure. We shall transfer the problem of dealing with
submodules of $M$ to $P$ where we know the claim is true by Lemma
\ref{T is b for fgp Umod}.

First, we shall show that $\cl_{\bnd}(f^{-1}(K)) =
f^{-1}(\cl_{\bnd}(K)).$

Let $x$ be in $\cl_{\bnd}(f^{-1}(K)).$ Then $g(x)=0,$ for every
$g\in\homo_{\U}(P, \U)$ such that $f^{-1}(K)\subseteq\ker g.$ We
need to show that $f(x)$ is in $\cl_{\bnd}(K),$ i.e. that
$h(f(x))=0$ for every $h\in\homo_{\U}(M, \U)$ with $K\subseteq\ker
h.$ Let $h$ be one such map. Letting $g = hf,$ we obtain a map in
$\homo_{\U}(P, \U)$ such that $g(f^{-1}(K))=hff^{-1}(K) = h(K)$
(since $f$ is onto). But $h(K)=0, $ and so $f^{-1}(K)\subseteq
\ker g.$ Hence, $g(x) = 0$ i.e. $h(f(x))=0.$

To show the converse, let $x$ be in $f^{-1}(\cl_{\bnd}(K)).$ Then
$h(f(x)) = 0$ for every $h$ $\in$ $\homo_{\U}(M, \U)$ such that
$K\subseteq\ker h.$ We need to show that $g(x) = 0$ for every
$g\in\homo_{\U}(P, \U)$ such that $f^{-1}(K)\subseteq\ker g.$ Let
$g$ be one such map. Since $f^{-1}(0)\subseteq
f^{-1}(K)\subseteq\ker g,$ we have $\ker f\subseteq\ker g.$ This
condition enables us to define a homomorphism $h: M\rightarrow \U$
such that $h(f(p)) = g(p)$ for every $p\in P.$ Then $h(K) = h(
f(f^{-1}(K))) = g(f^{-1}(K)) = 0,$ and so $h(f(x)) = 0.$ But this
gives us that $g(x) =0.$

It is easy to see that $f: P\rightarrow M$ induces an isomorphism
$P/f^{-1}(\cl_{\bnd}(K))$ $\rightarrow$ $M/\cl_{\bnd}(K).$ But $
\cl_{\bnd}(f^{-1}(K)) =f^{-1}(\cl_{\bnd}(K)),$ so we obtain that
$M/\cl_{\bnd}(K)$ is finitely generated projective (since
$P/\cl_{\bnd}(f^{-1}(K))$ is). So
$0\rightarrow\cl_{\bnd}(K)\rightarrow M\rightarrow
M/\cl_{\bnd}(K)\rightarrow 0 $ splits.

(2) To show that $\dim_{\U}(K) = \dim_{\U}(\cl_{\bnd}(K)),$ let us
look at a surjection $f:P\rightarrow M$ as in (1) and the
following two short exact sequences:
\[
\begin{array}{ccccccccc}
0 & \rightarrow & \ker f & \rightarrow & f^{-1}(K)& \rightarrow &
K
&  \rightarrow & 0, \nonumber \\
0 & \rightarrow & \ker f & \rightarrow & f^{-1}(\cl_{\bnd}(K)) &
\rightarrow & \cl_{\bnd}(K) & \rightarrow  & 0 \nonumber
\end{array}
\]

$\cl_{\bnd}(f^{-1}(K))=\cl_{\T}(f^{-1}(K))$ by Lemma \ref{T is b
for fgp Umod}. The dimension of $\cl_{\T}(f^{-1}(K))$ is the same
as the dimension of $f^{-1}(K)$ by the definition of the closure
and the theory $(\T, \bigP).$ Since
$\cl_{\bnd}(f^{-1}(K))=f^{-1}(\cl_{\bnd}(K)),$ the two exact
sequences give us $\dim_{\U}(K)= \dim_{\U}(\cl_{\bnd}(K)).$

(3) Let $K=0$ in (1) and (2).

(4) Since $\T\subseteq\bnd$ we have that $\T M\subseteq\bnd M.$
But, since $\dim_{\U}(\bnd M)=0$ and $\T M$ is the largest
submodule of $M$ with dimension zero, $\bnd M\subseteq\T M.$ The
rest of (4) follows from (3) and (1).
\end{proof}

Now we can prove the following.
\begin{thm}
For the ring $\U,$ \[(\T, \bigP) = \mbox{ Lambek torsion theory =
Goldie torsion theory }= (\bnd, \unb).\] \label{T=b for U}
\end{thm}
\begin{proof} Since we know that $(\T,
\bigP)\leq$ $\tau_L$ = $\tau_G$ = $(\bnd, \unb)$, it is sufficient
to show that $\bnd \subseteq \T$. Proposition \ref{b splits in U}
gives us that $\bnd M = \T M$ for every finitely generated $M$. To
finish the proof it suffices to show that every $\tau_L$-torsion
module is in $\T$. Let $M$ be $\tau_L$-torsion. Then all
submodules of $M$ are bounded (see Example (7) in Section
\ref{Examples}). So, all finitely generated submodules of $M$ are
bounded, and, hence in $\T$. Since the dimension of $M$ is the
supremum of the dimensions of its finitely generated submodules,
the dimension of $M$ is zero. Hence, $M$ is in $\T.$
\end{proof}

In contrast to the situation $(\T, \bigP)\lneqq(\bnd, \unb)$ for
the ring $\A,$ we have that $(\T, \bigP)=(\bnd, \unb)$ for the
ring $\U.$

The ring $\U$ is Ore because every von Neumann regular ring is
Ore, so the classical ring of quotients exists. Also, $\U$ is
semihereditary, so we have that $\U \subseteq Q_{\mathrm{cl}}(\U)
\subseteq Q_{\mathrm{max}}(\U) = E(\U).$ But $\U$ is
self-injective so $E(\U) = \U.$ Hence, $\U = Q_{\mathrm{cl}}(\U) =
Q_{\mathrm{max}}(\U) = E(\U).$ So, the classical torsion theory of
$\U$ is trivial. This indicates that there is only one nontrivial
torsion theory of the ring $\U$ of interest for us: $(\T, \bigP) =
\tau_L = \tau_G = (\bnd, \unb).$

This theory is neither trivial nor improper in general. Let $\NZ$
be the group von Neumann algebra of the group $\Zset$ and $\UZ$
its algebra of affiliated operators. Example 8.34 in \cite{Lu5}
gives us a flat $\NZ$-module $M$ with dimension zero. Since $M$ is
flat, $\UZ\otimes_{\NZ} M$ is nontrivial. Since $M$ has dimension
zero, $\dim_{\UZ}(\UZ\otimes_{\NZ} M)=0$, and so $\UZ\otimes_{\NZ}
M$ is in $\T.$ Thus, $(\T, \bigP)$ is not trivial for $\UZ.$

The theory $(\T, \bigP)$ is not improper whenever $\A$ (and hence
$\U(\A)$) is nontrivial since $\U(\A)$ is a torsion-free module
and, hence, not in $\T.$

If one is not interested in the $\smallt$-part of a module over a
finite von Neumann algebra $\A,$ one can work with $\U(\A)$
instead of $\A.$ For applications to topology, that means that we
can work with algebra of affiliated operators if we are not
interested in Novikov-Shubin invariants. See section 4 in
\cite{Re} for details about $L^2$-invariants via an algebra of
affiliated operators.

\section{{\bf Torsion Theories and Semisimplicity}}
\label{semisimple}

In this section, we shall see that the vanishing of certain
torsion theories is equivalent with the semisimplicity of $\U.$
First we need the following result.

\begin{lem}
Let $n$ be any positive integer. For every submodule $P$ of
$\U^n,$ \[\dim_{\U}(P) = \dim_{\A}(P\cap \A^n).\] \label{lema o
dimenziji ug ideala}
\end{lem}
\begin{proof}
If $P$ is finitely generated, then $P=\U\otimes_{\A}(P\cap\A^n)$
by Theorem \ref{moja K 0 theorem} and so \[\dim_{\U}(P)
=\dim_{\U}(\U\otimes_{\A}(P\cap\A^n))= \dim_{\A}(P\cap \A^n).\]

If $P$ is not finitely generated, write $P$ as a directed union of
its finitely generated submodules $P_i,$ $i\in I.$ Then
$P\cap\A^n$ is direct union of $P_i\cap\A^n,$ $i\in I.$ Thus, we
have
\[
\dim_{\U}(P)  =  \sup_{i\in I}\;\dim_{\U}(P_i)
 =  \sup_{i\in I}\;\dim_{\A}(P_i\cap\A^n)
 =  \dim_{\A}(P\cap \A^n).
\]
\end{proof}

Now we can prove the result about the equivalence of the vanishing
of certain torsion theories and the semisimplicity of $\U.$

\begin{thm}
The following are equivalent:
\begin{enumerate}
\item $\U$ is semisimple.
\item $(\T, \bigP)$ for $\U$ is trivial.
\item $(\T, \bigP)$ for $\A$ is equal to $(\smallt, \p).$
\item The $\T\p$-part of every $\A$-module is zero.
\item The $\T\p$-part of every cyclic $\A$-module is zero.
\end{enumerate}
\label{ug semisimple second set}
\end{thm}
\begin{proof} We shall show that
$(1)\Rightarrow (2)\Rightarrow (3)\Leftrightarrow  (4)\Rightarrow
(5)\Rightarrow (1).$

$(1)\Rightarrow (2).$ If $\U$ is semisimple, all $\U$-modules are
projective, and hence in $\bigP.$ So, $\T=0.$

$(2) \Rightarrow (3).$ Since $\smallt\subseteq\T$ for $\A,$ it is
sufficient to show that every module from $\T$ is in $\smallt.$ If
$M$ is in $\T$, then $\dim_{\U}(\U\otimes_{\A} M ) = \dim_{\A}(M)
= 0,$ so $\U\otimes_{\A} M = 0$ by assumption that there are no
nontrivial zero-dimensional $\U$-modules. But $\U\otimes_{\A} M=
0$ means that $M = \smallt M,$ so $M$ is in $\smallt$.

$(3) \Leftrightarrow (4).$ (3) is equivalent with (4) since $\T\p
M \cong \p\T M= \T M/\smallt M$.

$(4) \Rightarrow (5).$ is trivial.

$(5) \Rightarrow (1).$ To show that $\U$ is semisimple, it is
sufficient to show that every left ideal in $\U$ is a direct
summand. Let $I$ be a left ideal in $\U.$ Then $\cl_{\T}(I)$ is a
direct summand of $\U$ (by Proposition \ref{b splits in U}). We
shall show that $I$ is a direct summand by showing that $I =
\cl_{\T}(I).$

Since $\cl_{\T}(I)$ is a direct summand of $\U$, $\cl_{\T}(I)\cap
\A$ is a direct summand of $\A$ by Theorem \ref{moja K 0 theorem}.
Denote by $J$ the left ideal $I\cap \A$ and by $\overline{J}$ the
left ideal $\cl_{\T}(I)\cap \A.$ We shall show that $\overline{J}
= \cl_{\T}(J).$

Since $I\subseteq\cl_{\T}(I),$ we have $J \subseteq \overline{J}.$
$\overline{J}$ is $\T$-closed by Proposition 6.32 from \cite{Lam}
and the fact that $(\T, \bigP)$= Goldie torsion theory for $\A$.
Since $\cl_{\T}(J)$ is the smallest closed submodule containing
$J$ we have that $\cl_{\T}(J)\subseteq\overline{J}.$
$\overline{J}/\cl_{\T}(J)$ is contained in a finitely generated
projective module $\A/\cl_{\T}(J).$ So, $\overline{J}/\cl_{\T}(J)$
is a module in $\bigP.$ To show it is trivial, it is sufficient to
show that its dimension vanishes. This is the case since
\[
\begin{array}{rcll}
\dim_{\A}(\cl_{\T}(J)) & = & \dim_{\A}(J) &
(\mbox{Def. of }\T, \cl_{\T}\mbox{ and Prop. \ref{propofdim}})\\
 & = & \dim_{\A}(I\cap\A) & (\mbox{definition of }J)\\
 & = & \dim_{\U}(I) & (\mbox{by Lemma \ref{lema o dimenziji ug
 ideala}})\\
 & = & \dim_{\U}(\cl_{\T}(I)) &
(\mbox{Def. of }\T, \cl_{\T}\mbox{ and Prop. \ref{propofdim} for }\U)\\
 & = & \dim_{\A}(\cl_{\T}(I)\cap\A) &
 (\mbox{by Lemma \ref{lema o dimenziji ug ideala}})\\
 & = & \dim_{\A}(\overline{J}) & (\mbox{definition of }\overline{J})
\end{array}
\]
Thus, $\overline{J} = \cl_{\T}(J).$

By Proposition \ref{properties of closure}, $\cl_{\T}(J)/J =
\T(\A/J)$ and $\A/\cl_{\T}(J) = \bigP(\A/J)$. $\bigP(\A/J)$ is a
finitely generated projective module so the inclusion
$\T(\A/J)\hookrightarrow\A/J$ is split. So, $\cl_{\T}(J)/J =
\T(\A/J)$ is cyclic. Its $\T\p$-part is trivial by assumption, and
so \[0 =
\T(\A/J)/\smallt(\A/J)\cong\cl_{\T}(J)/\cl_{\smallt}(J).\] Thus,
$\cl_{\T}(J)/J = \cl_{\smallt}(J)/J =\smallt(\A/J)$ is in
$\smallt.$

Since $\cl_{\T}(J)/J$ is in $\smallt,$ $\U\otimes_{\A}\cl_{\T}(J)/
J = 0$ and, hence \[\U\otimes_{\A}\cl_{\T}(J) = \U\otimes_{\A}J.\]

Thus, \[
\begin{array}{rcll}
\cl_{\T}(I) & = & \U\otimes_{\A}(\cl_{\T}(I)\cap\A) &
(\mbox{by Theorem \ref{moja K 0 theorem}})\\
 & = & \U\otimes_{\A}\overline{J} &
 (\mbox{definition of }\overline{J}) \\
 & = & \U\otimes_{\A}\cl_{\T}(J) &
 (\mbox{since }\overline{J} = \cl_{\T}(J))\\
 & = & \U\otimes_{\A}J & (\mbox{by what we just showed})\\
 & = & \U\otimes_{\A}(I\cap\A) & (\mbox{definition of }J)\\
 & \subseteq & I & (I\mbox{ is a left ideal})
\end{array} \]

But, since $I$ is contained in $\cl_{\T}(I)$, we have that
$\cl_{\T}(I) = I.$ So, $I$ is a direct summand in $\U.$ Thus, $\U$
is semisimple. This finishes the proof.
\end{proof}

In view of the $\smallt-\T\p-\bigP$ filtration, the vanishing of
the $\T\p$-part of each $\A$-module is equivalent with $\U$ being
semisimple. The vanishing of the $\smallt$-part of every module is
equivalent with $\A=\U.$ Indeed, $\U\otimes_{\A}\U/\A=0$ so
$\U/\A$ is in $\smallt.$ Hence, if $\smallt=0,$ $\U=\A.$ The
converse is easy: if $\U=\A,$ then $\smallt M = \tor^{\A}_1(\U/A,
M)=0$ for every $\A$-module $M$.

In the case when a finite von Neumann algebra of interest is a
group von Neumann algebra $\vng,$ the algebra of affiliated
operators $\ug$ is semisimple if the group $G$ is finite. It is
easy to see that for finite group $G,$ $\ug = \vng = \Cset G$ and
$\Cset G$ is semisimple. The converse also holds: finite $G$ is
the {\em only} case when $\ug$ is semisimple. The proof of this
claim can be found in \cite{Lu5} (see the solution of the exercise
9.11). Thus, Theorem \ref{ug semisimple second set} asserts that
the $\T\p$-part is present for a large class of group von Neumann
algebras.

\end{document}